\documentclass[oneside]{amsart}
\usepackage[mathscr]{eucal}
\usepackage{amssymb, amsmath, array, amscd}
\usepackage[dvips]{graphicx}
\newtheorem{theorem}{Theorem}[section]

\newtheorem{proposition}{Proposition}[section]
\newtheorem{lemma}{Lemma}[section]

\newtheorem{remark}{Remark}[section]
\newtheorem{example}{Example}[section]

\begin{document}

\title[SINGULAR LEVI-FLAT HYPERSURFACES]
{SINGULAR LEVI-FLAT HYPERSURFACES \\ AND CODIMENSION ONE FOLIATIONS}

\author{Marco Brunella}

\address{Marco Brunella, IMB - CNRS UMR 5584, 9 Avenue Savary,
21078 Dijon, France}

\begin{abstract}
We study Levi-flat real analytic hypersurfaces with
singularities. We prove that the Levi foliation on the regular part
of the hypersurface can be holomorphically extended, in a suitable
sense, to neighbourhoods of singular points.
\end{abstract}

\maketitle

\vskip 2truecm

\section{Introduction}

Let $X$ be a complex manifold of dimension $n$. A smooth real
hypersurface $M\subset X$ is said to be {\bf Levi-flat} if the
codimension one distribution
$$T^{\mathbb C}M = TM \cap J(TM) \subset TM$$
is integrable, in Frobenius' sense. It follows that $M$ is smoothly
foliated by immersed complex manifolds of dimension $n-1$. We shall
denote by ${\mathcal F}_M$ such a {\bf Levi foliation}.

If $M$ is moreover real analytic, then its local structure is very
well understood: according to E. Cartan (see for instance \cite[\S
1.7]{BER}), around each $p\in M$ we can find local holomorphic
coordinates $z_1,...,z_n$ such that $M=\{\Im m z_1 =0\}$, and
consequently the leaves of ${\mathcal F}_M$ are $\{ z_1=c\}$,
$c\in{\mathbb R}$. In particular, the Levi foliation ${\mathcal
F}_M$ extends on a neighbourhood of $p$ to a codimension one
holomorphic foliation, with leaves $\{ z_1=c\}$, $c\in{\mathbb C}$.
It is then easy to see that these local extensions glue together,
giving a codimension one holomorphic foliation ${\mathcal F}$ on
some neighbourhood of $M$.

In this paper we are interested in finding a similar structure in
the case of {\it singular} Levi-flat hypersurfaces, a subject
sistematically initiated by Burns and Gong \cite{B-G}. Let us
firstly recall some terminology.

A closed subset $M\subset X$ is {\bf real analytic} if it is locally
defined by the vanishing of a (finite) collection of real analytic
functions. A real analytic subset $M$ is irreducible if it cannot be
expressed as $M=M_1\cup M_2$, with both $M_1$ and $M_2$ real
analytic and different from $M$. Any real analytic subset can be
decomposed (on relatively compact open subsets) into a finite
collection of irreducible components. An irreducible real analytic
subset $M$ has a well defined dimension $\dim_{\mathbb R}M$, and it
can be decomposed as a disjoint union $M=M_{reg}\cup M_{sing}$,
where: (i) $M_{reg}$ is nonempty and open in $M$, and it is formed
by those points of $M$ around which $M$ is a smooth real analytic
submanifold of $X$ of dimension $\dim_{\mathbb R}M$; (ii) $M_{sing}$
is a real analytic subset, all of whose irreducible components have
dimension strictly smaller that $\dim_{\mathbb R}M$. Remark,
however, that $M_{reg}$ may fail to be dense in $M$, as in the well
known Whitney umbrella.

When $\dim_{\mathbb R}M = \dim_{\mathbb R}X - 1 = 2n-1$, or more
generally each irreducible component of $M$ has dimension $2n-1$, we
shall call $M$ a {\bf real analytic hypersurface}. In that case, we
shall say that $M$ is {\bf Levi-flat} if $M_{reg}$ is Levi-flat in
the usual sense (Frobenius' integrability).

If $M\subset X$ is a Levi-flat real analytic hypersurface then we
have on $M_{reg}$ the Levi foliation ${\mathcal F}_{M_{reg}}$. A
natural question is: does this foliation extend holomorphically on
some neighbourhood of the {\it closure} $\overline{M_{reg}}$? Of
course, we need here to work, at least, with {\it singular}
codimension one holomorphic foliations; see for instance \cite{Cer}
for the basic dictionary. Let us see some examples.

\begin{example}\label{ex1} {\rm
Take a nonconstant holomorphic function $f:X\to{\mathbb C}$ and set
$M=\{ \Im m f=0\}$. Then $M$ is Levi-flat and $M_{sing}$ is the set
of critical points of $f$ lying on $M$. Leaves of the Levi foliation
on $M_{reg}$ are given by $\{ f=c\}$, $c\in{\mathbb R}$. Of course,
this Levi foliation can be extended to a singular holomorphic
foliation on the full $X$, generated by the kernel of $df$. More
complicated examples can be obtained by taking $M=\{ F\circ f=0\}$
for some real analytic function $F:{\mathbb C}\to{\mathbb R}$, or by
starting with a meromorphic $f$ \cite{B-G}.}
\end{example}

\begin{example}\label{ex2} {\rm
In $X={\mathbb C}^2$ with coordinates $z=x+iy$, $w=s+it$, consider
the irreducible real analytic hypersurface $M$ defined by
$$M = \{ t^2 = 4(y^2+s)y^2 \} .$$
We have $M_{sing} = \{ t=y=0\}$, a totally real plane in ${\mathbb
C}^2$. Note that $\overline{M_{reg}}\cap M_{sing} = \{ t=y=0, s\ge
0\}$, a proper subset of $M_{sing}$. This hypersurface is Levi-flat:
on $M_{reg}$, the leaves of ${\mathcal F}_{M_{reg}}$ are the complex
curves
$$L_c = \{ w=(z+c)^2, \Im m z\not= 0 \}, \qquad c\in{\mathbb R}.$$
For every $c\in{\mathbb R}$ the closure $\overline{L_c} =
\{w=(z+c)^2\}$ cuts $M_{sing}$ along the real curve $\{
s=(x+c)^2\}$. Thus, each point of $M_{sing}\cap\{ s>0\}$ belongs to
two closures $\overline{L_{c_1}}$ and $\overline{L_{c_2}}$, whereas
each point of $M_{sing}\cap\{ s<0\}$ belongs to no one. The
foliation ${\mathcal F}_{M_{reg}}$ extends holomorphically to a
neighbourhood of $M_{reg}$: just choose a square root of $w$ on
${\mathbb C}^2\setminus\{ t=0, s\ge 0\}$ (which is a neighbourhood
of $M_{reg}$), and then take the holomorphic foliation generated
there by $dw = 2\sqrt{w}dz$. However, it is clear that ${\mathcal
F}_{M_{reg}}$ cannot be extended to neighbourhoods of points of
$\overline{M_{reg}}\cap M_{sing}$, due to some ``ramification''
along $M_{sing}$ (or, more precisely, along the real line $\{
y=t=s=0\}\subset M_{sing}$: the other points are not seriously
problematic, for around them $M$ is the union of two smooth
Levi-flat hypersurfaces intersecting transversely). In spite of
this, we could say that ${\mathcal F}_{M_{reg}}$ can be extended in
some weak or multiform sense, as solution of the implicit
differential equation $(\frac{dw}{dz})^2 = 4w$.}
\end{example}

Our main result is an extrapolation of this last example. The
philosophy which is behind is that, even if sometimes it is
impossible to extend ${\mathcal F}_{M_{reg}}$, it is always possible
to extend it in a ``microlocal'' sense, i.e. after lifting to the
cotangent bundle. In other words, ${\mathcal F}_{M_{reg}}$ can
always be extended as a (transcendental) implicit differential
equation, or web. This will be better explained in Section 2 below.
For the moment, we state our result in the form of a ``resolution
theorem'' for singular Levi-flat hypersurfaces.

\begin{theorem}\label{thm1}
Let $X$ be a complex manifold of dimension $n$, and let $M\subset X$
be a Levi-flat real analytic hypersurface. Then there exist a
complex manifold $Y$ of dimension $n$, a Levi-flat real analytic
hypersurface $N\subset Y$, a (singular) codimension one holomorphic
foliation ${\mathcal F}$ on $Y$ extending ${\mathcal F}_{N_{reg}}$,
and a holomorphic map $\pi : Y\to X$, such that for some open
$N_0\subset N$:
\begin{enumerate}
\item[(i)] $\pi\vert_{N_0} : N_0\to M_{reg}$ is an isomorphism;
\item[(ii)] $\pi\vert_{\overline{N_0}} : \overline{N_0}\to
\overline{M_{reg}}$ is a proper map.
\end{enumerate}
\end{theorem}

In the Example \ref{ex2} above, we may take $Y={\mathbb C}^2$, $N$ a
real hyperplane, and $\pi$ a quadratic map which sends complex lines
in $N$ to the curves $\overline{L_c}$, $c\in{\mathbb R}$. Remark
that here $N_0$ is not equal to $N_{reg}$, it is only a (open and
dense) part of it. Unfortunately, during this procedure we loss the
isolated part of $M_{sing}$. Remark that $\overline{M_{reg}}$, as
well as $\overline{N_0}$, is not a real analytic subset, it is only
{\it sub}analytic. It could be interesting to generalise the above
Theorem, in some way, to the subanalytic setting, but our method is
strictly inside the real analytic world.

The Theorem above can be seen as a first step toward a resolution
procedure for singularities of Levi-flat hypersurfaces: the map
$\pi$ is a sort of blow-up of $\overline{M_{reg}}\setminus M_{reg}$,
allowing afterwards the extension of the Levi foliation. The second
step, which for the moment requires $n\le 3$, consists in applying
the Desingularisation Theorem of Seidenberg ($n=2$) or Cano ($n=3$)
\cite{Cer} \cite{Can} to the foliation ${\mathcal F}$. The third,
and easy, step is the analysis of Levi-flat real analytic
hypersurfaces which are tangent to a simple singularity of
codimension one foliation. We leave the details to the interested
reader.

The proof of Theorem \ref{thm1} is based on quite elementary
considerations concerning the complexification(s) of real analytic
subsets (not only hypersurfaces) in complex manifolds, which in some
sense stem from the seminal work of Diederich and Fornaess
\cite{D-F} and which are also exploited in \cite{B-G}.

\section{Lifting to the cotangent bundle}

In this Section we give a more precise statement of our main result,
and we reduce its proof to a general statement on the intrinsic
complexification of certain real analytic subsets.

Consider, as before, a Levi-flat real analytic hypersurface $M$ in a
complex manifold $X$, $\dim_{\mathbb C}X=n$. Let $PT^*X$ be the
projectivised cotangent bundle of $X$: a ${\mathbb C}P^{n-1}$-bundle
over $X$, whose fibre $PT_x^*X$ over $x\in X$ will be thought as the
set of complex hyperplanes in $T_xX$. Denote by $\pi$ the projection
$PT^*X\to X$.

The regular part $M_{reg}$ of $M$ can be lifted to $PT^*X$: just
take, for every $x\in M_{reg}$, the complex hyperplane
$$T_x^{\mathbb C}M_{reg} = T_xM_{reg}\cap J(T_xM_{reg})\subset
T_xX.$$ Call $$M_{reg}'\subset PT^*X$$ this lifting of $M_{reg}$.
Remark that it is no more a hypersurface: its (real) dimension
$2n-1$ is half of the real dimension of $PT^*X$. However, it is
still ``Levi-flat'', in a sense which will be precised below.

Take now a point $y$ in the closure $\overline{M_{reg}'}$,
projecting on $X$ to a point $x\in\overline{M_{reg}}$.

\begin{lemma}\label{lemma1}
There exists, in a neighbourhood $U_y\subset PT^*X$ of $y$, a real
analytic subset $N_y$ of dimension $2n-1$ containing $M_{reg}'\cap
U_y$.
\end{lemma}

(The notation suggests that $U_y$ and $N_y$ should be more properly
considered as {\it germs} at $y$. Similarly for the Proposition
below).

\begin{proof}
We choose local coordinates $z_1,...,z_n,w_1,...,w_{n-1}$ around $y$
such that $z_1,...,z_n$ are coordinates around $x$ and
$w_1,...,w_{n-1}$ corresponds to the hyperplane $dz_n +
\sum_{j=1}^{n-1} w_jdz_j = 0$. Take a real analytic equation $f$ of
$M$ around $x$, with $df\not= 0$ on $M_{reg}$. Then the real
analytic subset defined around $y$ by the equations
$$f=0 \quad , \qquad w_j = \frac{\partial f}{\partial z_j}\big/
\frac{\partial f}{\partial z_n} \quad , \quad j=1,...,n-1$$ contains
$M_{reg}'$ as an open subset. Indeed, $T^{\mathbb C}M_{reg}$ is the
Kernel of $\partial f$. This analytic subset may have dimension
larger than $2n-1$ (for instance, it contains the full fibres over
$M_{sing}$), but using a stratification of it we see that it
contains a real analytic subset $N_y$ (a stratum) of dimension
$2n-1$ and which still contains $M_{reg}'$.
\end{proof}

We choose such a $N_y$ as the minimal one (as a germ at $y$), by
taking intersections. We have $M_{reg}'\cap U_y\subset (N_y)_{reg}$,
but the inclusion may be strict, as in the Example \ref{ex2} of the
Introduction. Also, $N_y$ may have several irreducible components,
but each one contains a part of $M_{reg}'\cap U_y$, by minimality.

The crucial fact is now the following.

\begin{proposition}\label{prop1}
There exists, in a neighbourhood $V_y\subset U_y$ of $y$, a complex
analytic subset $Y_y$ of (complex) dimension $n$ containing $N_y\cap
V_y$.
\end{proposition}

Let us deduce Theorem \ref{thm1} from this Proposition.

First of all, $N_y\cap V_y$ is a real analytic hypersurface in
$Y_y$, and it is Levi-flat because each irreducible component
contains a Levi-flat piece \cite[Lemma 2.2]{B-G}. On $Y_y$, we have
a natural codimension one holomorphic (singular) distribution of
hyperplanes, given by the restriction of the canonical contact
structure of $PT^*X$. At points of $M_{reg}'$, this codimension one
distribution coincides (by tautology of the contact form) with
$T^{\mathbb C}M_{reg}$. In other words and more precisely, around
points of $M_{reg}'\cap V_y$ the projection $\pi :Y_y\to X$ is a
biholomorphism, sending $M_{reg}'$ to $M_{reg}$, thus $Y_y$ appears
as the graph of a local holomorphic section of $PT^*X$ extending the
real analytic section $M_{reg}\ni x' \mapsto y'=T_{x'}^{\mathbb
C}M_{reg}\in M_{reg}'$. This local holomorphic section is a local
hyperplane distribution on $X$. When lifted to $Y_y$ by the same
$\pi$, this distribution becomes a distribution on $Y_y$ coinciding
with the restriction of the contact structure.

Now, this distribution is integrable, because it is integrable on a
real hypersurface, hence it defines a codimension one (singular)
foliation ${\mathcal F}_y$ on $Y_y$. Clearly, this foliation
${\mathcal F}_y$ extends the Levi foliation of $(N_y)_{reg}\cap
V_y$, because so is for $M_{reg}'\cap V_y$ and each irreducible
component of $N_y\cap V_y$ contains a portion of $M_{reg}'\cap V_y$.

These local constructions are sufficiently canonical to be patched
together, when $y$ varies on $\overline{M_{reg}'}$: if
$Y_{y_1}\subset V_{y_1}$ and $Y_{y_2}\subset V_{y_2}$ are as above,
with $M_{reg}'\cap (V_{y_1}\cap V_{y_2})\not= \emptyset$, then
$Y_{y_1}\cap (V_{y_1}\cap V_{y_2})$ and $Y_{y_2}\cap (V_{y_1}\cap
V_{y_2})$ have some common irreducible components containing
$M_{reg}'\cap (V_{y_1}\cap V_{y_2})$, and so $Y_{y_1}$ and $Y_{y_2}$
can be glued by identifying those components. In this way, we obtain
an analytic space $Y$ of dimension $n$, with a Levi-flat
hypersurface $N$ and a holomorphic foliation ${\mathcal F}$,
enjoying all the properties stated in Theorem \ref{thm1}; the map
$\pi : Y\to X$ is deduced from the projection of $PT^*X$ to $X$.
Finally, to get a smooth $Y$ we just replace the possibly singular
$Y$ with a resolution of it, over which $N$ and ${\mathcal F}$ can
be lifted.

Let us reformulate this proof from a slightly different point of
view.

On some neighbourhood $X_0 \subset X$ of $M_{reg}$ we have a
holomorphic foliation ${\mathcal F}_0$ extending the Levi foliation
${\mathcal F}_{M_{reg}}$. The graph of this foliation is a complex
manifold $Y_0\subset PT^*X$ which projects bihlomorphically to $X_0$
by $\pi$. This $Y_0$ contains $M_{reg}'$, as a Levi-flat
hypersurface. The main point is then to extend $Y_0$ to a
neighbourhood of $\overline{M_{reg}'}$, and this is the content of
Proposition \ref{prop1}. In this way, an extension problem for
foliations (or, more appropriately, for webs) has been transformed
into an extension problem for complex analytic subspaces.

\begin{remark} {\rm
The analytic subset $Y_y$ can be seen as an implicit differential
equation on $X$, around $x$. But there are two quite different
situations. If $\pi :Y_y\to X$ is a finite map, then (by
Weierstrass) $Y_y$ has an equation which is {\it polynomial} in the
vertical variables $w_1,...,w_{n-1}$, so that the associated
implicit differential equation is ``algebraic in the derivatives''.
In some sense, even if there is a ramification we could say that
this situation is ``not-too-singular'', from the differential
equation point of view. If on the contrary $\pi :Y_y\to X$ is not
finite over $x$, we don't know if such an algebraicity in the
derivatives persists, i.e. if $Y_y$ can be analytically continued to
a neighbourhood of the full fibre over $x$.}
\end{remark}

\begin{example}\label{ex3} {\rm
Consider the subset $M\subset {\mathbb C}^2$ defined by
$$M = \bigcup_{\lambda\in{\mathbb S}^1} \{ z_2=\lambda z_1 + \varphi
(\lambda ) z_1^2 \}$$ where $\varphi :{\mathbb S}^1\to{\mathbb C}$
is a real analytic function. Outside the origin (and close to it)
$M$ is real analytic, smooth, and Levi-flat. Thus $M\setminus\{ 0\}$
can be lifted to $PT^*{\mathbb C}^2 = {\mathbb C}^2\times{\mathbb
C}P^1$, and the closure of this lifting is a smooth real analytic
threefold $N$ which cuts the fibre over $0$ along a smooth circle.
The cooking recipe to obtain the complex surface $Y$ containing $N$
is the following. We differentiate the complex curves contained in
$M$,
$$dz_2 = (\lambda + 2\varphi (\lambda )z_1)dz_1 ,$$
and then we ``eliminate'' $\lambda$ from the two equations $z_2 =
\lambda z_1 + \varphi (\lambda )z_1^2$ and $\frac{dz_2}{dz_1} =
\lambda +2\varphi (\lambda )z_1$, thus obtaining an analytic
relation between $z_1$, $z_2$ and $w_1=-\frac{dz_2}{dz_1}$. It is
quite clear that such an equation will be algebraic (in $w_1$) if
and only if $\varphi$ is algebraic, i.e. $\varphi (\lambda )=
\sum_{j=-\ell}^{\ell} a_j \lambda^j$. However, this algebraicity
condition on $\varphi$ seems to be also necessary and sufficient for
the real analyticity of $M$ at the origin. Thus, if $M$ is real
analytic at $0$ then we obtain an implicit differential equation
which is algebraic in the derivatives, and the question of the
previous Remark remains unanswered. Note, however, that whatever
$\varphi$ is, the subset $N$ is everywhere analytic; hence our
method, which is based on the analyticity of $N$ and not of $M$, is
unfortunately rather useless for these type of subtle problems.}
\end{example}

Proposition \ref{prop1} is a particular case of a more general fact,
for which it is convenient to introduce some more terminology.

Let $Z$ be a complex manifold of dimension $m$. Let $N\subset Z$ be
a real analytic subset of dimension $2n-1$. We shall say that $N$ is
a {\bf Levi-flat real analytic subset} if its regular part is
foliated by complex submanifolds of (complex) dimension $n-1$. As in
the hypersurface case, we shall denote by ${\mathcal F}_{N_{reg}}$
this {\bf Levi foliation} on $N_{reg}$; it is real analytic and of
real codimension one (in $N_{reg}$). This definition can be
reformulated in the language of CR geometry \cite[Ch. I]{BER}:
$N_{reg}$ is a CR submanifold of CR dimension $n-1$ (or CR
codimension one) whose complex tangent bundle $T^{\mathbb C}N_{reg}
= TN_{reg}\cap J(TN_{reg})$ is integrable. It is easy to see, as in
the hypersurface case \cite[Lemma 2.2]{B-G}, that it is sufficient
to check the integrability only on some open nonempty subset
$N_0\subset N_{reg}$.

As in the hypersurface case, the local structure of $N_{reg}$ is
well understood \cite[\S 1.8]{BER}: around every $z\in N_{reg}$
there are local holomorphic coordinates $z_1,...,z_m$ on $Z$ such
that
$$N_{reg} = \{ \Im m z_1=0, z_{m-n}=...=z_m=0\} .$$
In particular, around $z$ there is a canonically defined complex
submanifold of $Z$, of dimension $n$, which contains $N_{reg}$:
$$Y_0 = \{ z_{m-n}=...=z_m=0 \} .$$
This $Y_0$ is called {\bf intrinsic complexification} of $N_{reg}$
(around $z$). It is the smallest complex submanifold containing
$N_{reg}$. Remark that $N_{reg}$ is a (smooth) Levi-flat {\it
hypersurface} in $Y_0$, and the Levi foliation ${\mathcal
F}_{N_{reg}}$ can be canonically extended to a holomorphic
codimension one foliation ${\mathcal F}_{Y_0}$ on $Y_0$. This is the
situation that we have already met in the proof of Theorem
\ref{thm1}.

Now the basic result, from which Proposition \ref{prop1} follows,
states that this intrinsic complexification can be analytically
prolonged around points of the closure $\overline{N_{reg}}$:

\begin{theorem}\label{thm2}
Let $Z$ be a complex manifold of dimension $m$, and let $N\subset Z$
be a Levi-flat real analytic subset of dimension $2n-1$. Then, for
every $z\in \overline{N_{reg}}$, there exists a neighbourhood
$V\subset Z$ of $z$ and a complex analytic subset $Y\subset V$ of
dimension $n$ which contains $N_{reg}\cap V$
\end{theorem}

The proof will be given in the next section, and it is based on the
following idea. There exists a second type of complexification of
$N$, which could be called {\bf extrinsic complexification}
\cite{Car}. It lives in a larger space, not inside $Z$, but it has
the advantage that it can be defined around every point of $N$, not
only $N_{reg}$. The intrinsic complexification appears as a
``projection'' of the extrinsic one, and from the analytic
extendability to singular points of the latter we shall deduce the
analytic extendability of the former. In a vague sense, this is
related to \cite[Appendix]{D-F}.

It is however important to realize that this result holds because
the CR codimension of $N_{reg}$ is one, {\it and no more}. Here is
an example showing that an analogous statement for the higher CR
codimensional case may fail. This is related to the fact that,
whereas a holomorphic map from a complex curve is always ``locally
proper'', the same is no more true for maps from higher dimensional
complex manifolds, due to the phenomenon of ``contraction of
divisors''.

\begin{example}\label{ex4} {\rm
We work in ${\mathbb C}^3$, with coordinates $z_1,z_2,z_3$. Let
$S\subset {\mathbb C}^3$ be the real analytic surface defined by
$$z_2= (\Re e z_1)z_1$$
$$z_3= e^{\Re e z_1}z_1$$
Note that $S$ is smooth, being the graph over the $z_1$-axis of a
smooth function $F:{\mathbb C}_{z_1}\to {\mathbb C}_{z_2,z_3}$. At
the point $0\in S$ we have $T_0^{\mathbb C}S=T_0S$, a complex line,
whereas at any other point $p\in S$, $p\not= 0$, we have
$T_p^{\mathbb C}S=\{ 0\}$. Indeed, $0\in{\mathbb C}_{z_1}$ is the
only point where the differential of $F$ is ${\mathbb C}$-linear.
Thus, $S_0 = S\setminus \{ 0\}$ is a CR submanifold of CR dimension
0 and CR codimension 2. The ``Levi foliation'' is here the foliation
by points.

Along $S_0$ we have the intrinsic complexification $Y_0\supset S_0$,
which is a complex surface. However, this complex surface cannot be
extended around the point $0$. To see this, observe that over $\{
z_1\not= 0\}$ the surface $Y_0$ has equation
$$z_3 = e^{z_2/z_1}z_1$$
which has an essential singularity at $z_1=0$. More geometrically,
if we blow up ${\mathbb C}^3$ at the origin, $\widetilde{{\mathbb
C}^3} \buildrel b\over\to {\mathbb C}^3$, then the closure of
$b^{-1}(Y_0)$ has a trace on the exceptional divisor
$b^{-1}(0)\simeq {\mathbb C}P^2$ which contains the complex curve
with equation $w_3=e^{w_2}$, in affine coordinates $w_2=z_2/z_1$,
$w_3=z_3/z_1$. The transcendency of this curve implies that $Y_0$
cannot be holomorphically extended at $0$.}
\end{example}

\section{Complexification of Levi-flat subsets}

We start the proof of Theorem \ref{thm2} by recalling some facts
about (extrinsic) complexification of real analytic subsets
\cite{Car} \cite[Ch. X]{BER}.

Let $Z$ be a complex manifold of dimension $m$. We shall denote by
$Z^*$ the complex manifold {\it conjugate} to $Z$, that is obtained
by replacing the complex structure $J$ of $Z$ with the opposite
complex structure $-J$. If $A\subset Z$ is a complex analytic
subset, then the {\it same} subset is also a complex analytic subset
of $Z^*$, but with the opposite complex structure; it will be denote
by $A^*$. The diagonal $\Delta\subset Z\times Z^*$ is a totally real
submanifold. The projections of $Z\times Z^*$ onto $Z$ and $Z^*$
will be respectively denoted by $\Pi$ and $\Pi^*$, and the {\it
antiholomorphic} involution $Z\times Z^*\to Z\times Z^*$,
$(z_1,z_2)\mapsto (z_2,z_1)$, will be denoted by $\jmath$.

Let $N\subset Z$ be a real analytic subset of dimension $k$, and fix
a point $z_0\in\overline{N_{reg}}$. Without loss of generality for
our purposes, we will suppose that the germ of $N$ at $z_0$ is
irreducible. Set
$$N^\Delta = \{ (z,z)\in Z\times Z^* \ \vert\ z\in N\} .$$
Then in some sufficiently small neighbourhood of $(z_0,z_0)$, {\it
which we may assume equal to $Z\times Z^*$ by restricting $Z$},
there exists an irreducible complex analytic subset $N^{\mathbb C}$
of dimension $k$ such that
$$N^{\mathbb C}\cap\Delta = N^\Delta .$$
This $N^{\mathbb C}$ is called {\bf complexification} of $N$ at
$z_0$. It is obtained by complexifying the real analytic equations
defining $N$, and it can be characterized as the smallest (germ at
$(z_0,z_0)$ of) complex analytic subset containing $N^\Delta$. It
enjoys the following fundamental symmetry property: $\jmath$ leaves
$N^{\mathbb C}$ invariant, and $N^\Delta$ is the fixed point set of
$\jmath\vert_{N^{\mathbb C}}$.

A word of caution. A real analytic subset is not necessarily
coherent \cite{Car}. This means that if $z_1\in\overline{N_{reg}}$
is another point, close to $z_0$, then we can still define the
complexification $\widetilde{N^{\mathbb C}}$ of $N$ at $z_1$ (if the
germ of $N$ at $z_1$ has several irreducible components, we
complexify each one), however it may happen that
$\widetilde{N^{\mathbb C}}$, defined in some neighbourhood of
$(z_1,z_1)$, is {\it smaller} than the restriction of $N^{\mathbb
C}$ to that neighbourhood. More precisely, it may happen that the
germ of $\widetilde{N^{\mathbb C}}$ at $(z_1,z_1)$ is not the full
germ of $N^{\mathbb C}$ at $(z_1,z_1)$, but only a union of certain
irreducible components of it. The situation is even worst if we take
$z_1\in N\setminus\overline{N_{reg}}$: the local dimension of $N$ at
$z_1$ is smaller than $k$, and thus the complexification of $N$ at
$z_1$ has also smaller dimension; the germ of $\widetilde{N^{\mathbb
C}}$ at $(z_1,z_1)$ is then a proper subset of the germ of
$N^{\mathbb C}$.

Suppose now that $N$ is Levi-flat, $k=2n-1$. Take a regular point
$z\in N_{reg}$. As we saw in the previous section, in some
neighbourhood $V_z\subset Z$ of $z$ there exists a unique smooth
complex submanifold $Y_z$ of dimension $n$, which contains
$N_{reg}\cap V_z$ and over which the Levi foliation extends as a
holomorphic codimension one foliation ${\mathcal F}_z$. For a good
choice of $V_z$, we may assume that the space of leaves $Y_z\big/
{\mathcal F}_z$ is a disc ${\mathbb D}$, in which the Levi leaves
correspond to points of ${\mathbb I}=(-1,1)\subset{\mathbb D}$.
Hence, we have a well defined Schwarz reflection
$$s_z : Y_z\big/ {\mathcal F}_z \to Y_z\big/ {\mathcal F}_z$$
with respect to ${\mathbb I}$.

Provided that $V_z$ is sufficiently small, the complexification of
$N\cap V_z$ at $z$ is a smooth complex submanifold $N_z^{\mathbb C}$
of $U_z=V_z\times V_z^*\subset Z\times Z^*$, of dimension $2n-1$,
which in fact coincides with an irreducible component of $N^{\mathbb
C} \cap U_z$. The following Lemma gives the precise relation between
$Y_z$ and $N_z^{\mathbb C}$.

\begin{lemma} \label{lemma2}
The projection $\Pi$ induces a smooth holomorphic fibration of
$N_z^{\mathbb C}$ over $Y_z$, and the fibre of $\Pi : N_z^{\mathbb
C}\to Y_z$ over $z'\in Y_z$ is the Schwarz reflection of the leaf of
${\mathcal F}_z$ through $z'$ (with the opposite complex structure).
\end{lemma}

\begin{proof}
In suitable local coordinates $z_1,...,z_m$ we have
$$N_{reg} = \{ \Im m z_1=0, z_{m-n}=...=z_m=0 \}$$
$$Y_z = \{ z_{m-n}=...=z_m=0 \} .$$
If $z'\in Y_z$, then the leaf through $z'$ is
$$L_{z'} = \{ z_1=a, z_{m-n}=...=z_m=0 \}$$
where $a$ is the $z_1$-coordinate of $z'$, and its Schwarz
reflection is
$$s_z(L_{z'})= \{ z_1=\bar a , z_{m-n}=...=z_m=0 \} .$$
Local coordinates on $Z\times Z^*$ are provided by $z_1,   ,z_m,\bar
z_1,...,\bar z_m$, and in these coordinates
$$N_z^{\mathbb C} = \{ z_1=\bar z_1, z_{m-n}=...=z_m=\bar
z_{m-n}=...=\bar z_m =0 \} .$$ We therefore see that $N_z^{\mathbb
C}$ projects by $\Pi$ to $Y_z$, and the fibre over $z'$ is
$$\{ \bar z_1 =a, \bar z_{m-n}=...=\bar z_m=0 \} = s_z(L_{z'})^* .$$
\end{proof}

Of course, a similar statement holds also for the second projection
$\Pi^*$. Thus, $N_z^{\mathbb C}$ is a smooth complex hypersurface in
$Y_z\times Y_z^*\subset Z\times Z^*$, with a double fibration over
$Y_z$ and $Y_z^*$.

Note that for every $z'\in Y_z$ we can also consider the fibre over
$z'$ of the full $N^{\mathbb C}$, i.e. the intersection $N^{\mathbb
C} \cap (\{ z'\}\times Z^*)$. This is an analytic subset of $Z^*$
which extends the leaf $s_z(L_{z'})^*\subset V_z^*$. In particular,
each leaf of ${\mathcal F}_z$ has an {\it analytic continuation} to
the full $Z$. This fundamental fact was discovered in \cite{D-F}, in
a different context, and it was our starting point.

Let us now see what happens at the singular point $z_0$.

For every $p\in N^{\mathbb C}$ consider the vertical fibre through
$p$
$$F_p^* = \{ q\in N^{\mathbb C}\ \vert\ \Pi (q)=\Pi (p)\}$$
and define
$$d(p) = \dim_p (F_p^*)$$
where $\dim_p$ denotes the local dimension at $p$ (if $F_p^*$ is
reducible at $p$, we take the maximal dimension among irreducible
components). According to Cartan or Remmert \cite[\S V.3]{Loj}, the
function $d$ on $N^{\mathbb C}$ is Zariski semicontinuous, thus
equal to some constant $d_0$ over some Zariski open and dense
$\breve N^{\mathbb C}\subset N^{\mathbb C}$ and strictly greater
than $d_0$ over $N^{\mathbb C}\setminus \breve N^{\mathbb C}$. By
Lemma \ref{lemma2}, we have $d_0=n-1$, and moreover $\breve
N^{\mathbb C}$ contains generic points of $N_{reg}^\Delta$ (here we
say ``generic'', and not ``all'', because at special points
$(z,z)\in N_{reg}^\Delta$ the germ of $N^{\mathbb C}$ could have a
second irreducible component different from $N_z^{\mathbb C}$).

Set $p_0=(z_0,z_0)$.

If $d(p_0)=d_0=n-1$, then by the Rank Theorem  of Remmert \cite[\S
V.6]{Loj} there exists a neighbourhood $U_{z_0}\subset Z\times Z^*$
of $p_0$ such that $\Pi (N^{\mathbb C}\cap U_{z_0})$ is a complex
analytic subset $Y_{z_0}$ of $V_{z_0} = \Pi (U_{z_0})$, of dimension
$n$. The situation is not very different from the one of Lemma
\ref{lemma2}. And we are just in the conclusion of Theorem
\ref{thm2}: obviously $Y_{z_0}$ contains the previously defined
$Y_z$, for every $z\in N_{reg}\cap V_{z_0}$.

Hence, we shall suppose from now on that $d(p_0)=l\ge n.$ Thus the
vertical fibre $F_{p_0}^*$ contains an irreducible component
$Y^*\subset Z^*$ of dimension $l$, passing through $p_0$. Its
conjugate $Y=\jmath (Y^*)\subset Z$ is an irreducible component of
the horizontal fibre $F_{p_0} = \{ q\in N^{\mathbb C}\ \vert\ \Pi^*
(q)=\Pi^* (p_0)\}$. In the following, we shall consider $Y$
sometimes as a (horizontal) subset of $N^{\mathbb C}$, sometimes as
a subset of the base $Z$.

Obviously, $Y\subset Z$ is contained in $\Pi (N^{\mathbb C})$.
Because $\dim N^{\mathbb C} = 2n-1$ and $d_0=n-1$, the image of
$N^{\mathbb C}$ by $\Pi$ is a countable union of local analytic
subsets of dimension $\le n$ (by iterated application of the Rank
Theorem). Thus $Y$ cannot have dimension larger than $n$, and so
$$\dim Y =n.$$

To complete the proof, we will show that $Y$ contains $\Pi
(N^{\mathbb C}\cap U_{z_0})$, for some neighbourhood $U_{z_0}$ of
$z_0$ (thus, after all, $Y$ is equal to that projection, and so $Y$
is in fact the unique irreducible component of dimension $n$ of the
horizontal fibre). Because $Y$ is closed, it is sufficient to verify
this inclusion for $\Pi (\breve N^{\mathbb C}\cap U_{z_0})$.

Note that $Y\subset N^{\mathbb C}$ is not contained in $N^{\mathbb
C}\setminus \breve N^{\mathbb C}$, for dimensional reasons: because
$Y$ is horizontal of dimension $n$ and each vertical fibre through
$N^{\mathbb C}\setminus \breve N^{\mathbb C}$ has dimension $\ge n$,
$N^{\mathbb C}$ would have dimension $\ge 2n$. Thus, $Y$ intersects
$\breve N^{\mathbb C}$. If $p=(z,z_0)$ is a point of $Y\cap \breve
N^{\mathbb C}$, then (by the Rank Theorem) the image by $\Pi$ of a
small neighbourhood of $p$ in $N^{\mathbb C}$ is an analytic subset
$Y_z\subset V_z\subset Z$ of dimension $n$. The germ of $Y_z$ at $z$
could have several irreducible components, but at least one of them
is contained in $Y$, because $Y\cap V_z\subset Y_z$ and $\dim Y =
\dim Y_z$. By a connectivity argument (recall that $N^{\mathbb C}$
is irreducible, and therefore $\breve N^{\mathbb C}$ is connected),
we see that $Y$ contains all the points arising from the projection
of $\breve N^{\mathbb C}$ to $Z$, as claimed.

Let us conclude by looking again at the (counter)example of the
previous Section.

\begin{example}\label{ex5} {\rm
Take again the real analytic surface $S$ in ${\mathbb C}^3$ with
equation
$$z_2 = (\Re e z_1)z_1$$
$$z_3 = e^{\Re e z_1}z_1$$
which, outside $0$, is a CR submanifold of CR dimension 0 and CR
codimension 2. Its complexification is the complex analytic surface
$S^{\mathbb C}$ in ${\mathbb C}^3\times {\mathbb C}^{3*}$ with
equation (setting $\bar z_j = w_j$)
$$z_2 = \frac{z_1+w_1}{2}z_1 \quad ,\qquad w_2 =
\frac{z_1+w_1}{2}w_1$$
$$z_3 = e^{\frac{z_1+w_1}{2}}z_1 \quad ,\qquad w_3 =
e^{\frac{z_1+w_1}{2}}w_1$$ The projection of $S^{\mathbb C}$ to the
first factor ${\mathbb C}^3$ is not analytic: something like
$z_3=e^{z_2/z_1}z_1$, for $z_1\not= 0$. The horizontal fibre through
$0\in S^{\mathbb C}$ is the complex curve $C\subset {\mathbb C}^3$
given by
$$z_2=\frac{1}{2}z_1^2 \quad ,\qquad z_3=e^{\frac{1}{2}z_1}z_1.$$
It is, of course, contained in the projection of $S^{\mathbb C}$ but
not equal to. The dimensional arguments used several times in the
proof of Theorem \ref{thm2} cannot be used here.}
\end{example}

\end{document}